\documentclass[11pt,english]{article}
\usepackage[T1]{fontenc}
\usepackage[latin9]{inputenc}
\usepackage{geometry}
\usepackage{color}
\usepackage{babel}
\usepackage{amsmath}
\usepackage{amsthm}
\usepackage{amssymb}
\usepackage[all]{xy}
\usepackage[unicode=true,pdfusetitle,
 bookmarks=true,bookmarksnumbered=false,bookmarksopen=false,
 breaklinks=false,pdfborder={0 0 0},pdfborderstyle={},backref=page,colorlinks=true]
 {hyperref}
\hypersetup{citecolor=blue, linkcolor=blue}

\makeatletter
\theoremstyle{plain}
\newtheorem{thm}{\protect\theoremname}[section]
\theoremstyle{remark}
\newtheorem{rem}[thm]{\protect\remarkname}
\theoremstyle{plain}
\newtheorem{prop}[thm]{\protect\propositionname}

\usepackage[leftcaption]{sidecap}

\sidecaptionvpos{figure}{c}

\usepackage{footmisc}

\newcommand{\FigBesBeg}[1][1.0]{%
 \let\MyFigure\figure
 \let\MyEndfigure\endfigure
 \renewenvironment{figure}[1]{\begin{SCfigure}[#1]##1}{\end{SCfigure}}}

\newcommand{\FigBesEnd}{%
 \let\figure\MyFigure
 \let\endfigure\MyEndfigure}

\makeatother

\providecommand{\propositionname}{Proposition}
\providecommand{\remarkname}{Remark}
\providecommand{\theoremname}{Theorem}

\begin{document}
\global\long\def\F{\mathbb{\mathbb{\mathbf{F}}}}%
 
\global\long\def\rk{\mathbb{\mathrm{rank}}}%
 
\global\long\def\crit{\mathbb{\mathrm{Crit}}}%
 
\global\long\def\eqdf{\stackrel{\mathrm{def}}{=}}%
 
\global\long\def\E{\overrightarrow{E}}%
 
\global\long\def\spec{\mathrm{Spec}}%
 
\global\long\def\fw{{\cal F}_{w}}%
 
\global\long\def\crit{\mathrm{Crit}}%
 
\global\long\def\cyr{\mathfrak{CR}}%
 
\title{A Note on the Trace Method for Random Regular Graphs}
\author{Joel Friedman~~~~Doron Puder}
\maketitle
\begin{abstract}
The main goal of this note is to illustrate the advantage of analyzing
the non-backtracking spectrum of a regular graph rather than the ordinary
spectrum. We show that by switching to non-backtracking spectrum,
the method of proof used in \cite{Puder2015} yields a bound of $2\sqrt{d-1}+\frac{2}{\sqrt{d-1}}$
instead of the original $2\sqrt{d-1}+1$ on the second largest eigenvalue
of a random $d$-regular graph.
\end{abstract}

\section{Introduction}

Let $\Gamma$ be a $d$-regular graph on $N$ vertices. The adjacency
matrix $A_{\Gamma}$ of $\Gamma$ has $N$ real eigenvalues
\[
d=\lambda_{1}\ge\lambda_{2}\ge\ldots\ge\lambda_{N}\ge-d.
\]
The first eigenvalue is the trivial eigenvalue $\lambda_{1}=d$ corresponding
to the constant eigenfunction. Denote by $\lambda\left(\Gamma\right)\eqdf\max\left(\lambda_{2},-\lambda_{N}\right)$
the maximal absolute value of a non-trivial eigenvalue. It is well
known that many properties of the graph can be measured by the value
of $\lambda\left(\Gamma\right)$. In particular, $\Gamma$ has better
expanding properties the smaller $\lambda\left(\Gamma\right)$ is
(see e.g., \cite{HLW06}). The Alon-Boppana bound states that $\lambda_{2}\left(\Gamma\right)\ge2\sqrt{d-1}-o_{N}\left(1\right)$
\cite{Alo86,Nil91}, and Alon conjectured that for random $d$-regular
graphs, $\lambda_{2}\left(\Gamma\right)$ is very close to $2\sqrt{d-1}$.
This conjecture was proven by the first named author \cite{Fri08}:
for every fixed $d\ge3$ and every $\varepsilon>0$, if $\Gamma$
is a uniformly random $d$-regular graph on $N$ vertices, then $\lambda\left(\Gamma\right)\le2\sqrt{d-1}+\varepsilon$
asymptotically almost surely (a.a.s.), namely, this holds with probability
tending to $1$ as $N\to\infty$. 

The proof in \cite{Fri08} uses the trace method, where the non-trivial
eigenvalues are bounded by counting closed cycles in the graph. (We
elaborate more in Section \ref{sec:The-proof-in-Pud15} below). This
method was also used in other works studying $\lambda\left(\Gamma\right)$
for random $\Gamma$: these include \cite{BroderShamir}, \cite{Fri91},
as well as the more recent, shorter proof of Alon's conjecture \cite{bordenave2019new}.
Another paper in this line of works is \cite{Puder2015} by the second
named author, which uses results from combinatorial group theory to
establish the weaker result that $\lambda\left(\Gamma\right)\le2\sqrt{d-1}+1$
a.a.s. 

The Hashimoto non-backtracking matrix $B_{\Gamma}$ of the graph $\Gamma$
is a matrix depicting the adjacency of oriented edges in $\Gamma$.
Famously, the Ihara-Bass formula relates the spectrum of $B_{\Gamma}$
with that of $A_{\Gamma}$ and shows that one determines the other:
see Section \ref{sec:The-Hashimoto-non-backtracking} for details.
Consequently, one can bound the largest non-trivial eigenvalue of
the non-backtracking matrix $B_{\Gamma}$ and deduce a bound on $\lambda\left(\Gamma\right)$.
This is indeed how the proof works in both proofs of Alon's conjecture
\cite{Fri08,bordenave2019new}.

However, \cite{Puder2015} bounds $\lambda\left(\Gamma\right)$ directly.
The main point of this note is that by passing to bounding the largest
non-trivial eigenvalue of the Hashimoto matrix $B_{\Gamma}$, the
method of proof in \cite{Puder2015} gives the following significantly
better bound:
\begin{thm}
\label{thm:main}Fix $d\ge3$, and let $\Gamma$ be a random $d$-regular
simple\footnote{A simple graph has no loops and parallel edges. In general, graphs
in this note are not assumed to be necessarily simple.} graph on $n$ vertices chosen at uniform distribution. Then a.a.s.,
\[
\lambda\left(\Gamma\right)\le2\sqrt{d-1}+\frac{2}{\sqrt{d-1}}.
\]
\end{thm}

This note should thus be thought of as an addendum to \cite{Puder2015}.
The original result in \cite{Puder2015} and its improvement in Theorem
\ref{thm:main} fall short of proving the full strength of Alon's
conjecture as in \cite{Fri08,bordenave2019new}. However, we find
this addendum interesting for two main reasons. First, it illustrates
what seems to be a fundamental advantage of analyzing the non-backtracking
spectrum rather than the ordinary one. Second, the method of proof
in \cite{Puder2015} and here is very different than the one in \cite{Fri08,bordenave2019new},
and may be used in other directions. For example, this method, including
the improvement suggested here, is also applied in \cite{hanany2020word},
which gives bounds on the second eigenvalue of random Schreier graphs
of the symmetric group. We remark that Section 1.2.3 of \cite{friedman2014relativized}
also mentions a similar improvement -- working with Hashimoto matrices
as opposed to adjacency matrices -- for the results of \cite{BroderShamir}
and of \cite{Fri91}.\\

The paper is organized as follows. Section \ref{sec:The-Hashimoto-non-backtracking}
recalls the Hashimoto non-backtracking matrix and the relation between
the ordinary spectrum and the non-backtracking spectrum of a regular
graph. In Section \ref{sec:The-proof-in-Pud15} we briefly describe
the original proof from \cite{Puder2015}, and in Section \ref{sec:Proof-of-main-Theorem}
give the adjustments needed to establish Theorem \ref{thm:main}.
Finally, Section \ref{sec:Random-coverings} contains some remarks
regarding the more general model of random coverings of a fixed base
graph. 

\section{The Hashimoto non-backtracking matrix and the Ihara-Bass formula\label{sec:The-Hashimoto-non-backtracking}}

Let $\Gamma$ be an undirected $d$-regular graph, not necessarily
simple, on $N$ vertices and let $A=A_{\Gamma}$ be its $N\times N$
adjacency matrix. Denote by $\E$ the set of oriented edges of $\Gamma$,
namely, each edge of $\Gamma$ appears twice in this set, once with
every possible orientation, so $\left|\E\right|=Nd$. For $e\in\E$,
we denote by $\overline{e}$ the same edge with the reverse orientation,
and by $h\left(e\right)$ and $t\left(e\right)$ the head and tail
of $e$, respectively. The \emph{Hashimoto} or \emph{non-backtracking}
matrix $B=B_{\Gamma}$ is a $\left|\E\right|\times\left|\E\right|$
$0$-$1$ matrix with rows and columns indexed by the elements of
$\E$. The $e,f$ entry is defined by
\[
B_{e,f}=\begin{cases}
1 & \mathrm{if}~t\left(e\right)=h\left(f\right)~\mathrm{and}~f\ne\overline{e},\\
0 & \mathrm{otherwise.}
\end{cases}
\]
The Ihara-Bass formula states that 
\begin{equation}
\det\left(I_{\E}-Bx\right)=\left(1-x^{2}\right)^{N\left(d/2-1\right)}\det\left(I_{N}-Ax+\left(d-1\right)x^{2}I_{N}\right).\label{eq:Ihara-Bass}
\end{equation}
In fact, a similar formula holds more generally for arbitrary finite
graphs -- see, e.g., \cite{kotani20002,rangarajan2018combinatorial}
for more details and proofs. The Ihara-Bass formula shows that the
spectra of $A$ and $B$ completely determine each other: 
\begin{equation}
\spec\left(B\right)=\left\{ \pm1\right\} \cup\left\{ \mu\,\middle|\,\mu^{2}-\lambda\mu+\left(d-1\right)=0,~\lambda\in\spec\left(A\right)\right\} .\label{eq:dictionary A-B}
\end{equation}
In particular, every eigenvalue $\lambda\in\spec\left(A\right)$ corresponds
to two eigenvalues $\frac{\lambda\pm\sqrt{\lambda^{2}-4\left(d-1\right)}}{2}\in\spec\left(B\right)$.
The trivial eigenvalue $\lambda_{1}=d$ corresponds to $d-1,1\in\spec\left(B\right)$.
Every eigenvalue $\lambda\in\spec\left(A\right)$ with $\left|\lambda\right|\ge2\sqrt{d-1}$
gives rise to two \emph{real} eigenvalues in $\left[-\left(d-1\right),-1\right]\cup\left[1,d-1\right]$,
while every eigenvalue with $\left|\lambda\right|<2\sqrt{d-1}$ corresponds
to two non-real eigenvalues lying on the circle of radius $\sqrt{d-1}$
around $0$ in $\mathbb{C}$. For a nice treatment of this dictionary
between $\spec\left(A\right)$ and $\spec\left(B\right)$ in the regular
case, consult \cite[Section 3]{lubetzky2016cutoff}.

In particular, we think of $d-1$ as the trivial eigenvalue of $B$,
which, again, corresponds to the constant eigenfunction on oriented
edges. By ordering the multiset of eigenvalues of $B$ by their absolute
value, we obtain
\begin{equation}
d-1=\left|\mu_{1}\right|\ge\left|\mu_{2}\right|\ge\ldots\ge\left|\mu_{2N}\right|=\left|\mu_{2N+1}\right|=\ldots=\left|\mu_{dN}\right|=1.\label{eq:evalues of B}
\end{equation}
We let $\mu\left(\Gamma\right)\eqdf\left|\mu_{2}\right|$ denote the
largest absolute value of a non-trivial eigenvalue. If $N\ge2$ then
$\mu\left(\Gamma\right)\in\left[\sqrt{d-1},d-1\right]$. Notice that
if $\mu\left(\Gamma\right)>\sqrt{d-1}$, in which case $\mu_{2}$
is real, then
\begin{equation}
\lambda\left(\Gamma\right)=\mu\left(\Gamma\right)+\frac{d-1}{\mu\left(\Gamma\right)}.\label{eq:mu vs. lambda of G}
\end{equation}

\section{The proof in \cite{Puder2015}\label{sec:The-proof-in-Pud15}}

As explained in the introduction to \cite{Puder2015}, using well-known
contiguity results for different models of random regular graphs,
when $d=2k$ is even, it is enough to prove Theorem \ref{thm:main}
in the \emph{permutation model}\footnote{\label{fn:odd d}To deal with odd values of $d$, one can show that
the probabilistic bound holding for random $\left(d+1\right)$-regular
graphs also holds for random $d$-regular graphs -- see \cite[Claim 6.1]{Puder2015}.}. In this model, a random $d$-regular graph on $N$ vertices is generated
by sampling $k=\frac{d}{2}$ independent uniformly random permutations
$\sigma_{1},\ldots,\sigma_{k}$ in the symmetric group $S_{N}$ and
constructing the corresponding Schreier graph depicting the action
of $S_{N}$ on $\left[N\right]\eqdf\left\{ 1,\ldots,N\right\} $ with
respect to $\sigma_{1},\ldots,\sigma_{k}$. Namely, the $N$ vertices
of the graph are labeled $1,\ldots,N$, and for every $1\le i\le N$
and every $1\le j\le k$, one adds an edge $\left(i,\sigma_{j}\left(i\right)\right)$
to the graph. The resulting graph may contain loops and parallel edges.

In order to use the trace method, we bound the number of closed walks
in the random graph $\Gamma$ generated in the permutation model.
We may direct and label the edges of $\Gamma$ by labeling $\left(i,\sigma_{j}\left(i\right)\right)$
by $j$ and directing it as follows:
\[
\xymatrix{i\ar@{->}[r]^{j~~} & \sigma_{j}\left(i\right)}
.
\]
Then, every (closed) walk of length $t$ in $\Gamma$ corresponds
to some \emph{word }in\emph{ $\left\{ \sigma_{1}^{\pm},\ldots,\sigma_{k}^{\pm1}\right\} ^{t}$}.
For example, the left-to-right walk 
\[
\xymatrix{4\ar@{->}[r]^{5} & 3\ar@{<-}[r]^{4} & 7\ar@{->}[r]^{2} & 10\ar@{<-}[r]^{5} & 3\ar@{->}[r]^{2} & 4}
\]
corresponds to the word $\sigma_{5}\sigma_{4}^{-1}\sigma_{2}\sigma_{5}^{-1}\sigma_{2}$.
Moreover, the number of closed walks in $\Gamma$ corresponding to
a given word is exactly the number of fixed points of the permutation
in $S_{N}$ defined by the word\footnote{Here, for ease of description, we compose permutations from left to
right, although this is completely immaterial in the analysis.} in $\sigma_{1},\ldots,\sigma_{k}$.

The proof in \cite{Puder2015} heavily depends on deep results in
combinatorial group theory proven in \cite{Puder2014,PP15}. These
works consider random permutations sampled by fixed words in the free
group $\F_{k}$ with basis $X=\left\{ x_{1},\ldots,x_{k}\right\} $.
Given a word $w\in\F_{k}$, the corresponding random permutation is
$w\left(\sigma_{1},\ldots,\sigma_{k}\right)$ where, as before, $\sigma_{1},\ldots,\sigma_{k}$
are independent, uniformly random permutations in $S_{n}$. For example,
if $w=x_{5}x_{4}^{-1}x_{2}x_{5}^{-1}x_{2}\in\F_{5}$, the random permutation
is $w\left(\sigma_{1},\ldots,\sigma_{5}\right)=\sigma_{5}\sigma_{4}^{-1}\sigma_{2}\sigma_{5}^{-1}\sigma_{2}$. 

For a word $w\in\F_{k}$, denote by $\fw\left(N\right)$ the random
variable counting the number of fixed points of the random permutation
$w\left(\sigma_{1},\ldots,\sigma_{k}\right)\in S_{N}$. Assume that
$d=2k$ is even and that $\Gamma$ is a random $d$-regular graph
on $N$ vertices in the permutation model. The first step in the trace
method for bounding $\lambda\left(\Gamma\right)$ is the observation
that for any even $t\in\mathbb{Z}_{\ge1}$,
\begin{eqnarray}
\mathbb{E}\left[\lambda\left(\Gamma\right)^{t}\right] & \le & \mathbb{E}\left[\sum_{i=2}^{N}\lambda_{i}\left(\Gamma\right)^{t}\right]=\mathbb{E}\left[\mathrm{tr}\left(A_{\Gamma}^{~t}\right)\right]-d^{t}=\sum_{w\in\left(X\cup X^{-1}\right)^{t}}\left(\mathbb{E}\left[\fw\left(N\right)\right]-1\right).\label{eq:trace method - first step}
\end{eqnarray}
As shown in \cite{Puder2014,PP15}, the asymptotic behaviour of the
expectation $\mathbb{E}\left[\fw\left(N\right)\right]$ depends on
an algebraic invariant of $w$, called the primitivity rank. This
invariant, denoted $\pi\left(w\right)$, is equal to the smallest
rank of a subgroup $H\le\F_{k}$ which contains $w$ as an imprimitive
element, namely, such that $w$ does not belong to any free basis
of $H$. There are no such subgroups if and only if $w$ is primitive
in $\F_{k}$, and in this case we define $\pi\left(w\right)=\infty$.
The possible values of $\pi$ in $\F_{k}$ are $0,1,\ldots,k$ and
$\infty$. 

The set of subgroups $H$ of rank $\pi\left(w\right)$ containing
$w$ as an imprimitive element is denoted $\crit\left(w\right)$.
It turn outs that $\grave{\left|\crit\left(w\right)\right|<\infty}$
for every $w$ \cite[Section 4]{PP15}. Theorem 1.8 in \cite{PP15}
states that 
\[
\mathbb{E}\left[\fw\left(N\right)\right]=1+\frac{\left|\crit\left(w\right)\right|}{N^{\pi\left(w\right)-1}}+O\left(\frac{1}{N^{\pi\left(w\right)}}\right).
\]
The big-O term was later effectivised in \cite[Proposition 5.1]{Puder2015}
to yield the following version of the theorem. Here $\left|w\right|$
denotes the length of $w$, namely, the number of letters when written
in reduced form in the given basis $X$ of $\F_{k}$.
\begin{thm}
\label{thm:PP15}Let $w\in\F_{k}$ and assume that $\left|w\right|=t$.
Then for every $N>t^{2}$,
\[
\mathbb{E}\left[\fw\left(N\right)\right]\le1+\frac{1}{N^{\pi\left(w\right)-1}}\left(\left|\crit\left(w\right)\right|+\frac{t^{2+2\pi\left(w\right)}}{N-t^{2}}\right)\le1+\frac{\left|\crit\left(w\right)\right|}{N^{\pi\left(w\right)-1}}\left(1+\frac{t^{2+2\pi\left(w\right)}}{N-t^{2}}\right).
\]
\end{thm}

To count closed walks of length $t$ in $\Gamma$, we therefore care
to know how many words of length $t$ there are in $\F_{k}$ of every
given primitivity rank. This, incorporated with the number of critical
subgroups, is given by the following proposition:
\begin{thm}
\label{thm:counting words}\cite[Proposition 4.3 and Theorem 8.2]{Puder2015}
Let $k\ge2$ and $m\in\left\{ 1,\ldots,k\right\} $. Then
\begin{equation}
\limsup_{t\to\infty}\left[\sum_{w\in\F_{k}\colon~\left|w\right|=t~\&~\pi\left(w\right)=m}\left|\crit\left(w\right)\right|\right]^{1/t}=\max\left(\sqrt{2k-1},2m-1\right).\label{eq:couting reduced words}
\end{equation}
\end{thm}

In fact, for $m\ge2$, the $\limsup$ in the theorem can be replaced
by ordinary $\lim$, and for $m=1$ it can be replaced by an ordinary
$\lim$ on even values of $t$. The cases not covered by Theorem \ref{thm:counting words}
are $\pi\left(w\right)\in\left\{ 0,\infty\right\} $. But $\pi\left(w\right)=0$
if and only if $w=1$, and $\pi\left(w\right)=\infty$ if and only
if $w$ is primitive, if and only if $\left|\crit\left(w\right)\right|=0$,
if and only if $\mathbb{E}\left[\fw\left(N\right)\right]=1$ for all
$N$. At any rate, words $w$ with $\pi\left(w\right)=\infty$ do
not contribute to the summation (\ref{eq:original proof core}) below\footnote{For completeness, let us mention that for $k\ge3$, the number of
primitive words of length $t$ behaves like $\left(2k-3\right)^{t}$
-- see \cite{puder2014growth}.}. 

Now we reach a step in \cite{Puder2015} which was required because
the original proof directly analyzed the ordinary spectrum and counted
closed walks with possible backtracking in $\Gamma$. Note that the
word corresponding to a closed walk in $\Gamma$ is reduced if and
only if the walk has no backtracking. So for arbitrary walks of length
$t$, we need to consider not only reduced words of length $t$ but
any words in $\left(X\cup X^{-1}\right)^{t}$. Using (\ref{eq:couting reduced words})
and the ``extended cogrowth formula'' -- \cite[Theorem 4.4]{Puder2015}
-- we obtain:
\begin{thm}
\label{thm:counting non-reduced words}\cite[Cor. 4.5]{Puder2015}
Let $k\ge2$ and $m\in\left\{ 0,1,\ldots,k\right\} $. Then
\[
\limsup_{t\to\infty}\left[\sum_{w\in\left(X\cup X^{-1}\right)^{t}\colon~\pi\left(w\right)=m}\left|\crit\left(w\right)\right|\right]^{1/t}=\begin{cases}
2\sqrt{2k-1} & \mathrm{if}~2m-1\le\sqrt{2k-1}\\
\frac{2k-1}{2m-1}+2m-1 & \mathrm{if}~2m-1\ge\sqrt{2k-1}.
\end{cases}
\]
\end{thm}

\begin{rem}
\label{rem:cogrowth paper}The proof of the extended cogrowth formula
never appeared explicitly in print. The writing of ``Notes on the
cogrowth formula: the regular, biregular and irregular cases'', which
is mentioned in the reference list in \cite{Puder2015}, was never
completed. One reason is that the discussion leading to the current
note took place already in the beginning of 2016. This discussion
led to the realization that the extended cogrowth formula was immaterial
for the current method of proof. The second named author still stands
behinds the statement of Theorem 4.4 in \cite{Puder2015}. The ideas
in that proof are not too far from some of the existing proofs to
Grigorchuk's cogrowth formula. 
\end{rem}

The final computation in \cite{Puder2015} goes as follows. With assumptions
as above, it follows from (\ref{eq:trace method - first step}) that
for all $t<\sqrt{N}$ even, 
\begin{eqnarray}
\mathbb{E}\left[\lambda\left(\Gamma\right)^{t}\right] & = & \sum_{w\in\left(X\cup X^{-1}\right)^{t}}\left(\mathbb{E}\left[\fw\left(N\right)\right]-1\right)=\sum_{m=0}^{k}\sum_{w\in\left(X\cup X^{-1}\right)^{t}\colon\pi\left(w\right)=m}\left(\mathbb{E}\left[\fw\left(N\right)\right]-1\right)\label{eq:original proof core}\\
 & \stackrel{\mathrm{Thm}~\ref{thm:PP15}}{\le} & \sum_{m=0}^{k}\frac{1}{N^{m-1}}\sum_{w\in\left(X\cup X^{-1}\right)^{t}\colon\pi\left(w\right)=m}\left|\crit\left(w\right)\right|\left(1+\frac{t^{2+2\pi\left(w\right)}}{N-t^{2}}\right).\nonumber 
\end{eqnarray}
If $t\approx c\log N$ for some constant $c=c\left(d\right)$, taking
the $t$-th root of both side and then taking the limit as $N\to\infty$,
we may use Theorem \ref{thm:counting non-reduced words} to estimate
the right hand side and deduce that $\mathbb{E}\left[\lambda\left(\Gamma\right)^{t}\right]^{1/t}\le2\sqrt{d-1}+0.835$
-- for details see Section 6.1 in \cite{Puder2015}. Using Markov's
inequality gives that a.a.s.~$\lambda\left(\Gamma\right)\le2\sqrt{d-1}+0.84$.
Finally, as explained in Footnote \ref{fn:odd d}, this can be used
to get a slightly weaker bound when $d$ is odd, and all together,
for every $d\ge3$, a random $d$-regular graph on $N$ vertices satisfies
a.a.s.~$\lambda\left(\Gamma\right)\le2\sqrt{d-1}+1$.

\section{Proof of Theorem \ref{thm:main}\label{sec:Proof-of-main-Theorem}}

To establish the bound in Theorem \ref{thm:main}, we adapt the proof
from Section \ref{sec:The-proof-in-Pud15} to the non-backtracking
spectrum. The trace of the $t$-th power $B^{t}$ of the Hashimoto
matrix $B=B_{\Gamma}$ is equal to the number of \emph{cyclically
non--backtracking closed walks} of length $t$ in $\Gamma$. Each
such walk consists of a sequence of $t$ edges $e_{1},\ldots,e_{t}$
so that $t\left(e_{i}\right)=h\left(e_{\left(i+1\right)\pmod t}\right)$
and $e_{\left(i+1\right)\pmod t}\ne\overline{e_{i}}$ for $i=1,\ldots,t$.
When $d=2k$ is even we may use the permutation model as in Section
\ref{sec:The-proof-in-Pud15} to sample a random $d$-regular graph
on $N$ vertices. Every cyclically non-backtracking closed walk corresponds
to a fixed point of a cyclically reduced word in the permutations
$\sigma_{1},\ldots,\sigma_{k}$. The total number of cyclically non-backtracking
closed walks in $\Gamma$ is, therefore, 
\[
\sum_{w\in\cyr_{t}\left(\F_{k}\right)}\fw\left(N\right),
\]
where $\cyr_{t}\left(\F_{k}\right)$ denotes the set of cyclically
reduced words of length $t$ in $\F_{k}$. There is an exact formula
for the number of such words:
\begin{prop}
\cite[Prop. 17.2]{mann2011groups}\label{prop:Mann} The number of
cyclically reduced words of length $t$ in $\F_{k}$ is 
\[
\left|\cyr_{t}\left(\F_{k}\right)\right|=\left(2k-1\right)^{t}+k+\left(-1\right)^{t}\left(k-1\right).
\]
\end{prop}

The trace method in the non-backtracking case is based on the following
equality:
\[
\sum_{\mu\in\spec\left(B\right)}\mu^{t}=\mathrm{tr}\left(B^{t}\right)=\#\left\{ \mathrm{\mathrm{cyclically~non}-}\mathrm{backtracking~walks~of~length~}t\right\} .
\]
If $t$ is even, for every real eigenvalue $\mu$, the summand $\mu^{t}$
is positive. Since every non-real eigenvalue $\mu$ lies on $\left\{ z\in\mathbb{C}\colon\left|z\right|=\sqrt{d-1}\right\} $,
the summand $\mu^{t}$ in this case has real part at least $-\sqrt{d-1}^{t}$.
Recall also that the trivial eigenvalue is $d-1$ and that (at least)
$N\left(d-2\right)+1$ out of the $Nd$ eigenvalues are $\pm1$. Hence,
recalling the notation (\ref{eq:evalues of B}), for $t$ even we
have
\begin{eqnarray*}
\mathrm{tr}\left(B^{t}\right) & = & \left(d-1\right)^{t}+\sum_{i=2}^{2N-1}\mu_{i}\left(\Gamma\right)^{t}+N\left(d-2\right)+1,
\end{eqnarray*}
so
\begin{eqnarray}
\mathrm{Re}\left[\mu_{2}\left(\Gamma\right)^{t}\right] & = & \mathrm{tr}\left(B^{t}\right)-\left(d-1\right)^{t}-\sum_{i=3}^{2N-1}\mathrm{Re}\left[\mu_{i}\left(\Gamma\right)^{t}\right]-N\left(d-2\right)-1\nonumber \\
 & \le & \left[\sum_{w\in\cyr_{t}\left(\F_{k}\right)}\fw\left(N\right)\right]-\left(d-1\right)^{t}+2N\sqrt{d-1}^{t}-N\left(d-2\right)-1\nonumber \\
 & \stackrel{\mathrm{Prop.~\ref{prop:Mann}}}{=} & \left[\sum_{w\in\cyr_{t}\left(\F_{k}\right)}\left(\fw\left(N\right)-1\right)\right]+d-1+2N\sqrt{d-1}^{t}-N\left(d-2\right)-1\nonumber \\
 & \le & \left[\sum_{m=1}^{k}\sum_{w\in\cyr_{t}\left(\F_{k}\right)\colon\pi\left(w\right)=m}\left(\fw\left(N\right)-1\right)\right]+2N\sqrt{d-1}^{t}.\label{eq:step one}
\end{eqnarray}
Note that the last summation starts with $m=1$ and not with $m=0$
because only the trivial word $w=1$ has primitivity rank $\pi\left(w\right)=0$,
and the empty word is the only (cyclically) reduced word giving $1$.

Theorem \ref{thm:counting words} implies that for any $\varepsilon>0$
\begin{equation}
\sum_{w\in\cyr_{t}\left(\F_{k}\right)\colon\pi\left(w\right)=m}\left|\crit\left(w\right)\right|\le\sum_{w\in\F_{k}\colon\left|w\right|=t~\&~\pi\left(w\right)=m}\left|\crit\left(w\right)\right|\le\left[\max\left(\sqrt{2k-1},2m-1\right)+\varepsilon\right]^{t}\label{eq:bound on words}
\end{equation}
for every large enough $t$. Taking expectations on both sides of
(\ref{eq:step one}), we obtain that for every $\varepsilon>0$ and
large enough $t$,
\begin{eqnarray}
\mathbb{E}\left[\mathrm{Re}\left[\mu_{2}\left(\Gamma\right)^{t}\right]\right] & \stackrel{\mathrm{Thm.}~\ref{thm:PP15}}{\le} & 2N\sqrt{d-1}^{t}+\sum_{m=1}^{k}\frac{1}{N^{m-1}}\left(1+\frac{t^{2+2m}}{N-t^{2}}\right)\sum_{w\in\cyr_{t}\left(\F_{k}\right)~\colon~\pi\left(w\right)=m}\left|\crit\left(w\right)\right|\nonumber \\
 & \stackrel{\eqref{eq:bound on words}}{\le} & 2N\sqrt{d-1}^{t}+\left(1+\frac{t^{2+2k}}{N-t^{2}}\right)\sum_{m=1}^{k}\frac{1}{N^{m-1}}\left[\max\left(\sqrt{2k-1},2m-1\right)+\varepsilon\right]^{t}.\label{eq:step two}
\end{eqnarray}
We will soon take $t$ to be a function of $N$ so that as $N\to\infty$,
$N^{1/t}\to c$ for a constant $c$ specified below. Then for every
$\varepsilon>0$ and every large enough $N$, 
\[
\left(1+\frac{t^{2+2k}}{N-t^{2}}\right)\cdot2\left(k+1\right)\le\left(1+\varepsilon\right)^{t}.
\]
Because the right hand side of (\ref{eq:step two}) is at most $\left(k+1\right)$
times the maximal summand (among the $k+1$ summands), we get that
for every $\varepsilon>0$ and large enough $N$,
\begin{align}
 & \mathbb{E}\left[\mathrm{Re}\left[\mu_{2}\left(\Gamma\right)^{t}\right]\right]\le\nonumber \\
 & ~~~~\left[\left(1+\varepsilon\right)\cdot\max\left(\left\{ N^{1/t}\sqrt{d-1}\right\} \cup\left\{ \frac{2m-1}{N^{\left(m-1\right)/t}}\,\middle|\,2m-1\in\left[\sqrt{d-1},d-1\right]\right\} \right)\right]^{t},\label{eq:step three}
\end{align}
where we used the observation that if $2m-1<\sqrt{d-1}$ then the
term corresponding to $m$ in (\ref{eq:step two}) is $\frac{\sqrt{d-1}^{t}}{N^{\left(m-1\right)}}$,
and is thus strictly smaller than the first term $2N\sqrt{d-1}^{t}$.
A simple analysis yields that, at least for large values of $d$,
the optimal value of $t=t\left(N\right)$ is such that 
\[
N^{1/t}\to e^{\frac{2}{e\sqrt{d-1}}}
\]
as $N\to\infty$. With this value, whenever $2m-1\in\left[\sqrt{d-1},d-1\right]$,
write $m=\beta\sqrt{d-1}$ with $\beta>\frac{1}{2}$. Then 
\begin{eqnarray*}
\frac{2m-1}{N^{\left(m-1\right)/t}} & = & \frac{2\beta\sqrt{d-1}-1}{\left(N^{1/t}\right)^{\beta\sqrt{d-1}-1}}<N^{1/t}\sqrt{d-1}\cdot\frac{2\beta}{\left(N^{1/t}\right)^{\beta\sqrt{d-1}}}\\
 & \approx & N^{1/t}\sqrt{d-1}\cdot\frac{2\beta}{e^{2\beta/e}}\le N^{1/t}\sqrt{d-1},
\end{eqnarray*}
where the last inequality follows as $\frac{2\beta}{e^{2\beta/e}}\le1$
with equality if and only if $\beta=e/2$. Therefore, with this value
of $t$, we obtain from (\ref{eq:step three}) that for every $\varepsilon>0$,
\[
\mathbb{E}\left[\mathrm{Re}\left[\mu_{2}\left(\Gamma\right)^{t}\right]\right]\le\left[\left(1+\varepsilon\right)\cdot\sqrt{d-1}\cdot e^{\frac{2}{e\sqrt{d-1}}}\right]^{t}
\]
for every large enough $N$. Recall that if $\mu_{2}\left(\Gamma\right)$
is non-real, then it has absolute value $\sqrt{d-1}$, and so we always
have $\mathrm{Re}\left[\mu_{2}\left(\Gamma\right)^{t}\right]\ge-\sqrt{d-1}^{t}$
for $t$ even. Therefore, for $x=\frac{2}{e\sqrt{d-1}}$, 
\[
\mathrm{Prob}\left\{ \mu\left(\Gamma\right)>\left(1+2\varepsilon\right)\cdot e^{x}\sqrt{d-1}\right\} \cdot\left[\left(1+2\varepsilon\right)e^{x}\sqrt{d-1}\right]^{t}-\sqrt{d-1}^{t}\le\mathbb{E}\left[\mathrm{Re}\left[\mu_{2}\left(\Gamma\right)^{t}\right]\right]\le\left[\left(1+\varepsilon\right)e^{x}\sqrt{d-1}\right]^{t},
\]
which yields that for every $\varepsilon>0$
\[
\mathrm{Prob}\left\{ \mu\left(\Gamma\right)>\left(1+2\varepsilon\right)\cdot\sqrt{d-1}\cdot e^{\frac{2}{e\sqrt{d-1}}}\right\} \underset{N\to\infty}{\to}0.
\]
Finally, by (\ref{eq:mu vs. lambda of G}), when $\mu\left(\Gamma\right)>\sqrt{d-1}$,
we have that $\lambda\left(\Gamma\right)=\mu\left(\Gamma\right)+\frac{d-1}{\mu\left(\Gamma\right)}$,
and as $e^{x}+e^{-x}<2e^{x^{2}/2}$ for $x>0$, we conclude that
\[
\mathrm{Prob}\left\{ \lambda\left(\Gamma\right)>2\sqrt{d-1}\cdot e^{\frac{2}{e^{2}\left(d-1\right)}}\right\} \underset{N\to\infty}{\to}0.
\]
As $d$ is even so far, we have $d\ge4$ and $\frac{2}{e^{2}\left(d-1\right)}\le0.1$,
and for $y<0.1$, $e^{y}<1+1.1y$, and 
\[
2\sqrt{d-1}\cdot e^{\frac{2}{e^{2}\left(d-1\right)}}\le2\sqrt{d-1}\left(1+\frac{2.2}{e^{2}\left(d-1\right)}\right)<2\sqrt{d-1}+\frac{0.6}{\sqrt{d-1}}.
\]
Hence,
\begin{thm}
\label{thm:For d even}For $d$ \textbf{\emph{even}}, $d\ge4$, a
random $d$-regular graph $\Gamma$ on $N$ vertices satisfies 
\[
\lambda\left(\Gamma\right)\le2\sqrt{d-1}+\frac{0.6}{\sqrt{d-1}}
\]
a.a.s.~as $N\to\infty$.
\end{thm}

By \cite[Claim 6.1]{Puder2015}, an a.a.s.~bound on $\lambda\left(\Gamma\right)$
for a $\left(d+1\right)$-regular random graph also holds a.a.s.~for
a $d$-regular graph. We conclude that if $d\ge3$ is odd, then a
random $d$-regular graph $\Gamma$ on $N$ vertices satisfies
\[
\lambda\left(\Gamma\right)\le2\sqrt{d}+\frac{0.6}{\sqrt{d}}\le2\sqrt{d-1}+\frac{1}{\sqrt{d-1}}+\frac{0.6}{\sqrt{d}}<2\sqrt{d-1}+\frac{2}{\sqrt{d-1}}.
\]
This completes the proof of Theorem \ref{thm:main}.$\qed$
\begin{rem}
As explained in \cite[Section 6.2]{Puder2015}, for small values of
$d$, the optimal value of $t=t\left(N\right)$ in the proof above
is different and leads to better bounds.
\end{rem}

\section{Random coverings of a fixed graph\label{sec:Random-coverings}}

A random $2k$-regular graph on $N$ vertices in the permutation model
can be thought of as a random $N$-degree covering space of the bouquet
with one vertex and $k$ loops -- see \cite[Section 1]{Puder2015}.
More generally, \cite{Puder2015} deals with random coverings of an
arbitrary finite connected graph. In this case, extending Alon's conjecture,
the first named author conjectured in \cite{Fri03} that for every
$\varepsilon>0$, a random $N$-degree covering of a fixed graph $\Delta$
satisfies that a.a.s.~all \emph{new} eigenvalues of the covering
are at most $\rho+\varepsilon$ in absolute value, where $\rho$ is
the spectral radius of the universal cover of $\Delta$. In this more
general case, \cite{Puder2015} provided the best bounds at the time
it was written (and see the references therein for earlier bounds).
Slightly later, this conjecture was proven when the base graph is
regular in \cite{friedman2014relativized} (later split into a series
of papers starting with \cite{friedman2019relativized}) and in \cite{bordenave2019new}.
More recently, the conjecture was proven in full, namely, for arbitrary
finite base graph, in \cite{bordenave2019eigenvalues}.

We remark that the improvement suggested in the current note applies
more generally to random covers of a fixed regular graph. In this
case, \cite[Thm 1.5]{Puder2015} states that the largest new eigenvalue
of a random $N$-cover of a fixed $d$-regular graph is a.a.s.~less
than $2\sqrt{d-1}+0.84$. This can be improved to the same statement
as in Theorem \ref{thm:For d even}, namely, to $2\sqrt{d-1}+\frac{0.6}{\sqrt{d-1}}$. 

In the irregular case, there is no direct dictionary between the spectrum
of the ordinary spectrum and that of the non-backtracking spectrum.
However, \cite{bordenave2019eigenvalues} uses a variety of non-backtracking
operators to prove the conjecture about random coverings of arbitrary
graphs. In fact, they manage to prove something much stronger than
the original conjecture in \cite{Fri03}, and show the new spectrum
of a random covering of $\Delta$ is contained a.a.s.~in an $\varepsilon$-neighborhood
of the spectrum of the covering tree.

\section*{Acknowledgments}

We would like to thank Liam Hanany for beneficial discussions. The
computations that led to this note were carried out during the workshop
``Spectrum of Random Graphs'' held in January 2016 in CIRM, Luminy,
France, and we thank the organizers Charles Bordenave, Alice Guionnet
and B�lint Vir�g for this opportunity. J.F.~was partially supported
by an NSERC Grant. D.P.~was supported by the Israel Science Foundation:
ISF grant 1071/16. This project has received funding from the European
Research Council (ERC) under the European Union\textquoteright s Horizon
2020 research and innovation programme (grant agreement No 850956).

\bibliographystyle{alpha}
\bibliography{A_Note_on_the_Trace_Method}

\begin{thebibliography}{HLW06}

\bibitem[Alo86]{Alo86}
N.~Alon.
\newblock Eigenvalues and expanders.
\newblock {\em Combinatorica}, 6(2):83--96, 1986.

\bibitem[BC19]{bordenave2019eigenvalues}
C.~Bordenave and B.~Collins.
\newblock Eigenvalues of random lifts and polynomials of random permutation
  matrices.
\newblock {\em Annals of Mathematics}, 190(3):811--875, 2019.

\bibitem[Bor19]{bordenave2019new}
C.~Bordenave.
\newblock A new proof of {F}riedman's second eigenvalue {T}heorem and its
  extension to random lifts.
\newblock {\em Annales Scientifiques de l'Ecole Normale Supérieure}, 2019.

\bibitem[BS87]{BroderShamir}
A.~Broder and E.~Shamir.
\newblock On the second eigenvalue of random regular graphs.
\newblock In {\em The 28th Annual Symposium on Foundations of Computer
  Science}, pages 286--294, 1987.

\bibitem[FK14]{friedman2014relativized}
J.~Friedman and D.~Kohler.
\newblock The relativized second eigenvalue conjecture of {A}lon.
\newblock preprint arXiv:1403.3462, 2014.

\bibitem[FK19]{friedman2019relativized}
J.~Friedman and D.~Kohler.
\newblock On the relativized {A}lon second eigenvalue conjecture {I}: Main
  theorems, examples, and outline of proof.
\newblock preprint arXiv:1911.05688, 2019.

\bibitem[Fri91]{Fri91}
J.~Friedman.
\newblock On the second eigenvalue and random walks in random $d$-regular
  graphs.
\newblock {\em Combinatorica}, 11(4):331--362, 1991.

\bibitem[Fri03]{Fri03}
J.~Friedman.
\newblock Relative expanders or weakly relatively {R}amanujan graphs.
\newblock {\em Duke Mathematical Journal}, 118(1):19--35, 2003.

\bibitem[Fri08]{Fri08}
J.~Friedman.
\newblock {\em A proof of {A}lon's second eigenvalue conjecture and related
  problems}, volume 195 of {\em Memoirs of the AMS}.
\newblock AMS, september 2008.

\bibitem[HLW06]{HLW06}
S.~Hoory, N.~Linial, and A.~Wigderson.
\newblock Expander graphs and their applications.
\newblock {\em Bulletin of the American Mathematical Society}, 43(4):439--562,
  2006.

\bibitem[HP20]{hanany2020word}
Liam Hanany and Doron Puder.
\newblock Word measures on symmetric groups.
\newblock 2020.

\bibitem[KS00]{kotani20002}
M.~Kotani and T.~Sunada.
\newblock Zeta functions of finite graphs.
\newblock {\em Journal of Mathematical Sciences-University of Tokyo},
  7(1):7--25, 2000.

\bibitem[LP16]{lubetzky2016cutoff}
E.~Lubetzky and Y.~Peres.
\newblock Cutoff on all {R}amanujan graphs.
\newblock {\em Geometric and Functional Analysis}, 26(4):1190--1216, 2016.

\bibitem[Man11]{mann2011groups}
A.~Mann.
\newblock {\em How groups grow}, volume 395.
\newblock Cambridge university press, 2011.

\bibitem[Nil91]{Nil91}
A.~Nilli.
\newblock On the second eigenvalue of a graph.
\newblock {\em Discrete Mathematics}, 91(2):207--210, 1991.

\bibitem[PP15]{PP15}
D.~Puder and O.~Parzanchevski.
\newblock Measure preserving words are primitive.
\newblock {\em Journal of the American Mathematical Society}, 28(1):63--97,
  2015.

\bibitem[Pud14]{Puder2014}
D.~Puder.
\newblock Primitive words, free factors and measure preservation.
\newblock {\em Israel J. Math.}, 201(1):25--73, 2014.

\bibitem[Pud15]{Puder2015}
D.~Puder.
\newblock Expansion of random graphs: new proofs, new results.
\newblock {\em Inventiones Mathematicae}, 201(3):845--908, 2015.

\bibitem[PW14]{puder2014growth}
D.~Puder and C.~Wu.
\newblock Growth of primitive elements in free groups.
\newblock {\em Journal of the London Mathematical Society}, 90(1):89--104,
  2014.

\bibitem[Ran18]{rangarajan2018combinatorial}
B.~Rangarajan.
\newblock A combinatorial proof of {I}hara-{B}ass's formula for the {Z}eta
  function of regular graphs.
\newblock In {\em 37th IARCS Annual Conference on Foundations of Software
  Technology and Theoretical Computer Science (FSTTCS 2017)}. Schloss
  Dagstuhl-Leibniz-Zentrum fuer Informatik, 2018.

\end{thebibliography}

\noindent Joel Friedman, Department of Computer Science, University
of British Columbia, Vancouver, BC V6T 1Z4, CANADA\\
\texttt{jf@cs.ubc.ca} \\

\noindent Doron Puder, School of Mathematical Sciences, Tel Aviv University,
Tel Aviv, 6997801, Israel\\
\texttt{doronpuder@gmail.com}
\end{document}